\documentclass[12pt,reqno]{amsart}
\usepackage{amssymb}
\begin{document}
\theoremstyle{plain}
\newtheorem{thm}{Theorem}[section]
\newtheorem{theorem}[thm]{Theorem}
\newtheorem{lemma}[thm]{Lemma}
\newtheorem{corollary}[thm]{Corollary}
\newtheorem{corollary and definition}[thm]{Corollary and Definition}
\newtheorem{proposition}[thm]{Proposition}
\newtheorem{example}[thm]{Example}
\theoremstyle{definition}
\newtheorem{notation}[thm]{Notation}
\newtheorem{claim}[thm]{Claim}
\newtheorem{remark}[thm]{Remark}
\newtheorem{remarks}[thm]{Remarks}
\newtheorem{conjecture}[thm]{Conjecture}
\newtheorem{definition}[thm]{Definition}
\newtheorem{problem}[thm]{Problem}
\newcommand{\Diff}{{\rm Diff}}
\newcommand{\fz}{\frak{z}}
\newcommand{\zar}{{\rm zar}}
\newcommand{\an}{{\rm an}}
\newcommand{\red}{{\rm red}}
\newcommand{\codim}{{\rm codim}}
\newcommand{\rank}{{\rm rank}}
\newcommand{\Pic}{{\rm Pic}}
\newcommand{\Div}{{\rm Div}}
\newcommand{\Hom}{{\rm Hom}}
\newcommand{\im}{{\rm im}}
\newcommand{\Spec}{{\rm Spec}}
\newcommand{\sing}{{\rm sing}}
\newcommand{\reg}{{\rm reg}}
\newcommand{\Char}{{\rm char}}
\newcommand{\Tr}{{\rm Tr}}
\newcommand{\Gal}{{\rm Gal}}
\newcommand{\Min}{{\rm Min \ }}
\newcommand{\Max}{{\rm Max \ }}
\newcommand{\soplus}[1]{\stackrel{#1}{\oplus}}
\newcommand{\dlog}{{\rm dlog}\,}    
\newcommand{\limdir}[1]{{\displaystyle{\mathop{\rm
lim}_{\buildrel\longrightarrow\over{#1}}}}\,}
\newcommand{\liminv}[1]{{\displaystyle{\mathop{\rm
lim}_{\buildrel\longleftarrow\over{#1}}}}\,}
\newcommand{\boxtensor}{{\Box\kern-9.03pt\raise1.42pt\hbox{$\times$}}}
\newcommand{\sext}{\mbox{${\mathcal E}xt\,$}}
\newcommand{\shom}{\mbox{${\mathcal H}om\,$}}
\newcommand{\coker}{{\rm coker}\,}
\renewcommand{\iff}{\mbox{ $\Longleftrightarrow$ }}
\newcommand{\onto}{\mbox{$\,\>>>\hspace{-.5cm}\to\hspace{.15cm}$}}
\catcode`\@=11
\def\opn#1#2{\def#1{\mathop{\kern0pt\fam0#2}\nolimits}}
\def\bold#1{{\bf #1}}%
\def\underrightarrow{\mathpalette\underrightarrow@}
\def\underrightarrow@#1#2{\vtop{\ialign{$##$\cr
 \hfil#1#2\hfil\cr\noalign{\nointerlineskip}%
 #1{-}\mkern-6mu\cleaders\hbox{$#1\mkern-2mu{-}\mkern-2mu$}\hfill
 \mkern-6mu{\to}\cr}}}
\let\underarrow\underrightarrow
\def\underleftarrow{\mathpalette\underleftarrow@}
\def\underleftarrow@#1#2{\vtop{\ialign{$##$\cr
 \hfil#1#2\hfil\cr\noalign{\nointerlineskip}#1{\leftarrow}\mkern-6mu
 \cleaders\hbox{$#1\mkern-2mu{-}\mkern-2mu$}\hfill
 \mkern-6mu{-}\cr}}}
\let\amp@rs@nd@\relax
\newdimen\ex@
\ex@.2326ex
\newdimen\bigaw@
\newdimen\minaw@
\minaw@16.08739\ex@
\newdimen\minCDaw@
\minCDaw@2.5pc
\newif\ifCD@
\def\minCDarrowwidth#1{\minCDaw@#1}
\newenvironment{CD}{\@CD}{\@endCD}
\def\@CD{\def\A##1A##2A{\llap{$\vcenter{\hbox
 {$\scriptstyle##1$}}$}\Big\uparrow\rlap{$\vcenter{\hbox{%
$\scriptstyle##2$}}$}&&}%
\def\V##1V##2V{\llap{$\vcenter{\hbox
 {$\scriptstyle##1$}}$}\Big\downarrow\rlap{$\vcenter{\hbox{%
$\scriptstyle##2$}}$}&&}%
\def\={&\hskip.5em\mathrel
 {\vbox{\hrule width\minCDaw@\vskip3\ex@\hrule width
 \minCDaw@}}\hskip.5em&}%
\def\verteq{\Big\Vert&&}%
\def\noarr{&&}%
\def\vspace##1{\noalign{\vskip##1\relax}}\relax\let\amp@rs@nd@&\iffalse}\fi
 \CD@true\vcenter\bgroup\relax\let\\=\cr\iffalse}\fi\tabskip\z@skip\baselineskip20\ex@
 \lineskip3\ex@\lineskiplimit3\ex@\halign\bgroup
 &\hfill$\m@th##$\hfill\cr}
\def\@endCD{\cr\egroup\egroup}
\def\>#1>#2>{\amp@rs@nd@\setbox\z@\hbox{$\scriptstyle
 \;{#1}\;\;$}\setbox\@ne\hbox{$\scriptstyle\;{#2}\;\;$}\setbox\tw@
 \hbox{$#2$}\ifCD@
 \global\bigaw@\minCDaw@\else\global\bigaw@\minaw@\fi
 \ifdim\wd\z@>\bigaw@\global\bigaw@\wd\z@\fi
 \ifdim\wd\@ne>\bigaw@\global\bigaw@\wd\@ne\fi
 \ifCD@\hskip.5em\fi
 \ifdim\wd\tw@>\z@
 \mathrel{\mathop{\hbox to\bigaw@{\rightarrowfill}}\limits^{#1}_{#2}}\else
 \mathrel{\mathop{\hbox to\bigaw@{\rightarrowfill}}\limits^{#1}}\fi
 \ifCD@\hskip.5em\fi\amp@rs@nd@}
\def\<#1<#2<{\amp@rs@nd@\setbox\z@\hbox{$\scriptstyle
 \;\;{#1}\;$}\setbox\@ne\hbox{$\scriptstyle\;\;{#2}\;$}\setbox\tw@
 \hbox{$#2$}\ifCD@
 \global\bigaw@\minCDaw@\else\global\bigaw@\minaw@\fi
 \ifdim\wd\z@>\bigaw@\global\bigaw@\wd\z@\fi
 \ifdim\wd\@ne>\bigaw@\global\bigaw@\wd\@ne\fi
 \ifCD@\hskip.5em\fi
 \ifdim\wd\tw@>\z@
 \mathrel{\mathop{\hbox to\bigaw@{\leftarrowfill}}\limits^{#1}_{#2}}\else
 \mathrel{\mathop{\hbox to\bigaw@{\leftarrowfill}}\limits^{#1}}\fi
 \ifCD@\hskip.5em\fi\amp@rs@nd@}
\newenvironment{CDS}{\@CDS}{\@endCDS}
\def\@CDS{\def\A##1A##2A{\llap{$\vcenter{\hbox
 {$\scriptstyle##1$}}$}\Big\uparrow\rlap{$\vcenter{\hbox{%
$\scriptstyle##2$}}$}&}%
\def\V##1V##2V{\llap{$\vcenter{\hbox
 {$\scriptstyle##1$}}$}\Big\downarrow\rlap{$\vcenter{\hbox{%
$\scriptstyle##2$}}$}&}%
\def\={&\hskip.5em\mathrel
 {\vbox{\hrule width\minCDaw@\vskip3\ex@\hrule width
 \minCDaw@}}\hskip.5em&}
\def\verteq{\Big\Vert&}
\def\novarr{&}
\def\noharr{&&}
\def\SE##1E##2E{\slantedarrow(0,18)(4,-3){##1}{##2}&}
\def\SW##1W##2W{\slantedarrow(24,18)(-4,-3){##1}{##2}&}
\def\NE##1E##2E{\slantedarrow(0,0)(4,3){##1}{##2}&}
\def\NW##1W##2W{\slantedarrow(24,0)(-4,3){##1}{##2}&}
\def\slantedarrow(##1)(##2)##3##4{%
\thinlines\unitlength1pt\lower 6.5pt\hbox{\begin{picture}(24,18)%
\put(##1){\vector(##2){24}}%
\put(0,8){$\scriptstyle##3$}%
\put(20,8){$\scriptstyle##4$}%
\end{picture}}}
\def\vspace##1{\noalign{\vskip##1\relax}}\relax\let\amp@rs@nd@&\iffalse}\fi
 \CD@true\vcenter\bgroup\relax\let\\=\cr\iffalse}\fi\tabskip\z@skip\baselineskip20\ex@
 \lineskip3\ex@\lineskiplimit3\ex@\halign\bgroup
 &\hfill$\m@th##$\hfill\cr}
\def\@endCDS{\cr\egroup\egroup}
\newdimen\TriCDarrw@
\newif\ifTriV@
\newenvironment{TriCDV}{\@TriCDV}{\@endTriCD}
\newenvironment{TriCDA}{\@TriCDA}{\@endTriCD}
\def\@TriCDV{\TriV@true\def\TriCDpos@{6}\@TriCD}
\def\@TriCDA{\TriV@false\def\TriCDpos@{10}\@TriCD}
\def\@TriCD#1#2#3#4#5#6{%
\setbox0\hbox{$\ifTriV@#6\else#1\fi$}
\TriCDarrw@=\wd0 \advance\TriCDarrw@ 24pt
\advance\TriCDarrw@ -1em
\def\SE##1E##2E{\slantedarrow(0,18)(2,-3){##1}{##2}&}
\def\SW##1W##2W{\slantedarrow(12,18)(-2,-3){##1}{##2}&}
\def\NE##1E##2E{\slantedarrow(0,0)(2,3){##1}{##2}&}
\def\NW##1W##2W{\slantedarrow(12,0)(-2,3){##1}{##2}&}
\def\slantedarrow(##1)(##2)##3##4{\thinlines\unitlength1pt
\lower 6.5pt\hbox{\begin{picture}(12,18)%
\put(##1){\vector(##2){12}}%
\put(-4,\TriCDpos@){$\scriptstyle##3$}%
\put(12,\TriCDpos@){$\scriptstyle##4$}%
\end{picture}}}
\def\={\mathrel {\vbox{\hrule
   width\TriCDarrw@\vskip3\ex@\hrule width
   \TriCDarrw@}}}
\def\>##1>>{\setbox\z@\hbox{$\scriptstyle
 \;{##1}\;\;$}\global\bigaw@\TriCDarrw@
 \ifdim\wd\z@>\bigaw@\global\bigaw@\wd\z@\fi
 \hskip.5em
 \mathrel{\mathop{\hbox to \TriCDarrw@
{\rightarrowfill}}\limits^{##1}}
 \hskip.5em}
\def\<##1<<{\setbox\z@\hbox{$\scriptstyle
 \;{##1}\;\;$}\global\bigaw@\TriCDarrw@
 \ifdim\wd\z@>\bigaw@\global\bigaw@\wd\z@\fi
 \mathrel{\mathop{\hbox to\bigaw@{\leftarrowfill}}\limits^{##1}}
 }
 \CD@true\vcenter\bgroup\relax\let\\=\cr\iffalse}\fi
 \tabskip\z@skip\baselineskip20\ex@
 \lineskip3\ex@\lineskiplimit3\ex@
 \ifTriV@
 \halign\bgroup
 &\hfill$\m@th##$\hfill\cr
#1&\multispan3\hfill$#2$\hfill&#3\\
&#4&#5\\
&&#6\cr\egroup%
\else
 \halign\bgroup
 &\hfill$\m@th##$\hfill\cr
&&#1\\%
&#2&#3\\
#4&\multispan3\hfill$#5$\hfill&#6\cr\egroup
\fi}
\def\@endTriCD{\egroup}
\newcommand{\sA}{{\mathcal A}}
\newcommand{\sB}{{\mathcal B}}
\newcommand{\sC}{{\mathcal C}}
\newcommand{\sD}{{\mathcal D}}
\newcommand{\sE}{{\mathcal E}}
\newcommand{\sF}{{\mathcal F}}
\newcommand{\sG}{{\mathcal G}}
\newcommand{\sH}{{\mathcal H}}
\newcommand{\sI}{{\mathcal I}}
\newcommand{\sJ}{{\mathcal J}}
\newcommand{\sK}{{\mathcal K}}
\newcommand{\sL}{{\mathcal L}}
\newcommand{\sM}{{\mathcal M}}
\newcommand{\sN}{{\mathcal N}}
\newcommand{\sO}{{\mathcal O}}
\newcommand{\sP}{{\mathcal P}}
\newcommand{\sQ}{{\mathcal Q}}
\newcommand{\sR}{{\mathcal R}}
\newcommand{\sS}{{\mathcal S}}
\newcommand{\sT}{{\mathcal T}}
\newcommand{\sU}{{\mathcal U}}
\newcommand{\sV}{{\mathcal V}}
\newcommand{\sW}{{\mathcal W}}
\newcommand{\sX}{{\mathcal X}}
\newcommand{\sY}{{\mathcal Y}}
\newcommand{\sZ}{{\mathcal Z}}
\newcommand{\A}{{\mathbb A}}
\newcommand{\B}{{\mathbb B}}
\newcommand{\C}{{\mathbb C}}
\newcommand{\D}{{\mathbb D}}
\newcommand{\E}{{\mathbb E}}
\newcommand{\F}{{\mathbb F}}
\newcommand{\G}{{\mathbb G}}
\newcommand{\HH}{{\mathbb H}}
\newcommand{\I}{{\mathbb I}}
\newcommand{\J}{{\mathbb J}}
\newcommand{\M}{{\mathbb M}}
\newcommand{\N}{{\mathbb N}}
\renewcommand{\P}{{\mathbb P}}
\newcommand{\Q}{{\mathbb Q}}
\newcommand{\R}{{\mathbb R}}
\newcommand{\T}{{\mathbb T}}
\newcommand{\U}{{\mathbb U}}
\newcommand{\V}{{\mathbb V}}
\newcommand{\W}{{\mathbb W}}
\newcommand{\X}{{\mathbb X}}
\newcommand{\Y}{{\mathbb Y}}
\newcommand{\Z}{{\mathbb Z}}
\title{Improvement on Parameters of Algebraic-Geometry Codes from Hermitian Curves}   
\author{Siman Yang}
\thanks{The first author was partially supported by the program for Chang Jiang Scholars and Innovative Research Team in University}
\begin{abstract}
Motivated by Xing's method [7], we show that there exist $[n,k,d]$
linear Hermitian codes over $\F_{q^2}$ with $k+d\geq n-3$ for all
sufficiently large $q$. This improves the asymptotic bounds given
in $[9,10]$.
\end{abstract}

\maketitle


\ \\

\section{Introduction}
We first review Goppa's construction and properties of Algebraic-Geometry codes using the language of global function fields.
A global function field $F$ over a finite field $\F_q$ is a function field having the following properties:
\begin{enumerate}
\item[(i)] $F$ is an algebraic function field with constant field $\F_q$.
\item[(ii)] $\F_q$ is algebraically closed in $F$.
\end{enumerate}

A place of $F$ is an equivalence class of valuations of $F$. A place of $F$ is rational if it has degree $1$, i.e., its residue class field is $\F_q$. A divisor of $F$ is a formal sum $\sum_{P} n_P P$  with $n_P \in \Z$ and all but finitely many $n_P=0$. For a nonzero function $f\in F$, the principal divisor div$(f)$ is $\sum_{P} \nu _P(f) P$, where $\nu _P$ is the normalized discrete valuation corresponding to the place $P$.

Denoted by $N(F)$ the number of rational places of a function field $F/\F_q$ and $g(F)$ the genus of $F/\F_q$.
By the Hasse-Weil bound (cf. [4]) one has
$$
N(F) \leq q+1+2g\sqrt{q}. \eqno{(1.1)}
$$

Let $F/\F_q$ be a global function field of genus $g$ with $n$ rational places, and $G$ is a divisor of $F$.
The Riemann-Roch space associated to $G$ consists of sections of $G$, i.e., $L(G)=\{f\in F: \mbox{div}(f)+G\geq 0\}\cup\{0\}$ is a finite dimensional vector space over $\F_q$.
$P_1, \ldots, P_n$ are distinct rational places of $F$ with $\{P_1, \ldots, P_n\}\cap $supp$(G)=\emptyset$.
An Algebraic-Geometry code is defined as the image of the $\F_q$-linear map from $L(G)$ to $\F_q^n : f \mapsto (f(P_1),\ldots, f(P_n))$.
If $g\leq deg(G)<n$, one gets linear $[n,k,d]$ code over $\F_q$ with $k=\mbox{dim}(L(G))\geq \mbox{deg}(G)+1-g$ (by the Riemann-Roch theorem),  and $d\geq n-\mbox{deg}(G)$ (because a nonzero section of $G$ has at most deg$G$ zeros), cf. [6].
Thus, $k+d\geq n+1-g$.

Let $C$ be a linear code over an alphabet $\F_q$. Denoted by $n(C)$, $k(C)$, $d(C)$ the length, dimension, and minimum distance of $C$ respectively. We call the ratio $k(C)/n(C)$ the transmission rate and $d(C)/n(C)$ the error-detection rate of the code.  It is a theoretical important task in coding theory to construct a sequence of codes with length goes to infinity together with asymptotic positive transmission rate and error-detection rate. Goppa's construction of Algebraic-Geometry codes leads to the Tsfasman-Vladut-Zink bound [1] which is a breakthrough in coding theory as it beats the Gilbert-Varshamov bound in an open interval over a finite field of size $q_0^2$ with $q_0\geq 7$. This motivates the construction of codes with parameters better than Goppa's construction.

However, for about twenty years after Goppa's construction, the refinement mainly concerned exhibiting suitable curves with ``many'' rational points (with respect to genus) and algorithmic improvement to the resulting codes. Recently, Xing [8] gave a new construction of nonlinear algebraic-geometry codes by exploiting the sections' derivatives to find codes with better asymptotic parameters than Goppa's. Elkies [2,3] estimated the size of the set of rational sections of bounded degree of the line bundle $L_D$ associated to a degree zero divisor $D$ on a curve $C$ to construct algebraic-geometry codes. His constructions also led to an improved asymptotic bound better than the Tsfasman-Vladut-Zink bound. Very recently, Stichtenoth, Niederreiter, Ozbudak and some others also constructed several nonlinear codes with better parameters than Goppa's construction. In this paper we employ an idea of Xing [7] to yield linear codes from Hermitian curves with asymptotic parameters better than in the record (see [11], [9], [10]). Like all those codes mentioned above the codes we present are nonconstructive.

A maximal curve $F/\F_q$ of genus $g$ is a curve achieving the Hasse-Weil bound.
Clearly, $q$ must be a square if there exists a maximal curve of positive genus over $\F_q$. Hermitian curves $F/\F_{q^2}$ are a class of important examples of maximal curves defined by an equation of the form
$$
Y^q +Y = X^{q+1}. \eqno{(1.2)}
$$
It is known that $g(F)=(q^2-q)/2$ and $N(F)=q^3+1$, among them, there are $q^3$ affine $\F_{q^2}$-rational points $P_1, \ldots, P_{q^3}$ and one infinite $\F_{q^2}$-rational point $P_{\infty}$ on this curve. It is obvious that Hermitian curves are maximal curves.

Yang and Kumar [11] determined the true minimum distance of one-point Hermitian codes as follows:
For any integer $2g-1<t<q^3$, the one-point code $C_L(tP_{\infty}, D)$, where $D = P_1 +\ldots + P_{q^3}$, is a $[q^3, t-g+1, d]$ linear code over $\F_{q^2}$ with $d=q^3-t+c$ for some $0\leq c<q$. This improved the Goppa's estimate of mimimum distance of Hermitian codes to $O(q)$. By employing a method of Xing [7], Xing and Chen [9], Xu[10] improved the Goppa's estimate of mimimum distance of Hermitian codes to $O(q^2)$ in some range respectively.  The major effort here is expended on finding a specific divisor $G$ of a prescribed degree and satisfies $L(G-D') =\{0\}$ for all large subset $D'\subseteq \{P_1,\ldots , P_{q^3}\}$ of a prescribed degree. Our main result is as follows.

\begin{theorem}
There exist $[n,k,d]$ linear codes over $\F_{q^2}$ with $k+d\geq n-3$ for all sufficiently large $n=q^3$, where $q$ is a prime power.
\end{theorem}
This improves the minimum distance of sufficiently large length Hermitian code very close to the genus of the relied curve than Goppa's construction. We will use the following lemma in [7].

\section{Codes from Hermitian Curves}
\begin{lemma}
Let $F$ be an algebraic curve with $n$ rational points $P_1, \ldots, P_n$.
For fixed positive integers $s>m$ define
$$
N_{s,m}=|\{\sum_{P\in I} P +D : I\subseteq \{P_1,\ldots ,P_{n}\}, |I|=m, D\geq 0, \mbox{D}=s-m\}|. \eqno{(2.1)}
$$
Suppose $N_{s,m}<h(F)$ (denoted by $h(F)$ the class number of $F$). Then there exists a divisor $G$ of degree $s$ with $\{P_1, \ldots, P_n\}\cap $ supp$(G)=\emptyset$ and the AG code $C_L(P_1, \ldots, P_n;G)$ is $[n,k,d]$ code with $k=s-g+1$ and $d\geq n-m+1$.
\end{lemma}

Obviously, the above construction can improve $s-m+1$ of Goppa's estimate on minimum distance of Algebraic-geometry codes. We next search for such integers $s$ and $m$ satifying  Eq.(2.1) with large value $s-m$ for a Hermitian curve $F/\F_{q^2}$. To estimate $N_{s,m}$ one needs to count the number of positive divisors of fixed degree. Thus it is natural to consider the zeta-function of $F$ which is defined by
$$
Z_F(T)=\sum _{i=0}^{\infty} A_i T^i,
$$
where $A_i$ is defined as the number of positive divisors of $F$ of degree $i$ for all intergers $i\geq 0$. It is well known that the zeta-function of a curve over $\F_{q^2}$ is a rational function of the form
$$
Z_F(T)=\frac{L_F(T)}{(1-T)(1-q^2 T)}
$$
We have the following upper bound of $A_k$.
\begin{proposition}
For a maximal curve over $\F_{q^2}$ with genus $g$, it follows that
$$
A_k=\sum _{i=0}^k \binom{2g}{i}\frac{q^{2k+2-2i}-1}{q^2-1}q^i < h(F)q^{2k+2-2g}.
$$
\end{proposition}

\begin{proof}
The formula of $A_k$ follows directly from the formula of L-function of a maximal curve over $\F_{q^2}$, which is $L_F(T)=(1+qT)^{2g}$, cf. [5]. Note that $h(F)=L_F(1)$, thus one obtain the above inequality by some simple calculations.
\end{proof}

Combining Lemma 2.1, we yield the following result:
\begin{proposition}
Suppose $F$ is a maximal curve over $\F_{q^2}$ of genus $g$ with at least $n$ rational points and for some integers $l$ and $t$ holds
$$
\binom{n}{l} q^{2t+2-2g}\leq 1.
$$
Then there exists $[n,k,d]$ code over $\F_{q^2}$ with $k=l+t-g+1$ and $d\geq n-l+1$.
\end{proposition}

By applying Proposition 2.3 we are ready to prove our main theorem.

\begin{theorem}
There exist $[n,k,d]$ linear codes over $\F_{q^2}$ with $k+d\geq n-3$ for all sufficiently large $n=q^3$, where $q$ is a prime power.
\end{theorem}

\begin{proof}
Below we construct Hermitian codes over $\F_{q^2}$ of sufficiently large length $n=q^3$. Fix $n$, we take $l=\lfloor n\alpha \rfloor$, $t=g-1+\lfloor (\theta -1)l \rfloor$ for some $0<\theta , \alpha <1$. We compute
\begin{align*}
&\displaystyle\lim _{n\rightarrow \infty} \log _q ^{\binom{n}{l} q^{2t+2-2g}}/n\\
=&-\alpha \log _q ^{\alpha} -(1-\alpha)\log _q ^{1-\alpha} +2\alpha (\theta -1).
\end{align*}

To achieve that the above sum is negative, one must choose $\theta < 1- \log _q ^2 \frac{H_2(\alpha)}{2\alpha}$,
where $H_2(\delta)=-\delta \log _2 ^{\delta}-(1-\delta) \log _2 ^{1-\delta}$ is the $2$-th entropy function.
Taking $\theta = 1- \log ^2 _q \frac{H_2(\alpha)}{2\alpha}-o(1)$ and $\alpha =q^{\varepsilon -3}$ for some positive $\varepsilon=\varepsilon _q $
with $\displaystyle\lim _{q\rightarrow \infty   }q^{\varepsilon}\rightarrow 2$.
By Proposition 2.3 we achieve the following inequalities

\begin{align*}
d_{\text{improved}}:= d- d_{\text{Goppa}}&\geq g-1-\frac{\log ^2 _q H_2(\alpha)}{2\alpha}l-o(l)\\
&=g-1+\frac{q^{\varepsilon}}{2}(\varepsilon -3)+\frac{q^3-q^{\varepsilon}}{2}\log _q ^{1-q^{\varepsilon -3}}-o(q^{\varepsilon})\\
&=g-1+\frac{q^{\varepsilon}}{2}(\varepsilon -3)+\frac{-2q^{\varepsilon}+q^{-3+2\epsilon}+\ldots}{4\ln q}-o(2)\\
&=g-4-o(1).
\end{align*}

\end{proof}

\begin{remark}
We set $\displaystyle\lim _{q\rightarrow \infty   }q^{\varepsilon}\rightarrow 2$ such that the codes we yield are not trivial codes. One can show that this value can be replaced by any value larger than one.
\end{remark}

\begin{remark}
In [10, Remark 4.6], the improvement to Goppa's construction for Hermitian codes is close to $g-q$. Clearly, our construction yields codes with better parameters over large fields. We note that most results in this article work for any maximal curve.

\end{remark}

\vskip .2in
\noindent
Siman Yang\\
Department of Mathematics, East China Normal University,\\
500, Dongchuan Rd., Shanghai, P.R.China 200241. \ \ e-mail:
smyang@math.ecnu.edu.cn
 \\ \\


\begin{thebibliography}{D}

\bibitem[1] {}M. A. Tsfasman, S.G. Vladut, T. Zink, {\it Modular curves, Shimura curves, and Goppa codes, better than the Varshamov-Gilbert bound\/}. Math. Nachr. {\bf 109} (1982), 21--28.


\bibitem[2] {}N.D. Elkies, {\it Excellent nonlinear codes from modular curves\/}. In STOC 01: Proceedings of the 33rd Annual ACM Symposium on Theory of Computing, (2001), 200--208.

\bibitem[3] {}N.D. Elkies, {\it Still better nonlinear codes from modular curves\/}. http://arxiv.org/abs/math/0308046.

\bibitem[4] {}H. Stichtenoth, {\it Algebraic Function Fields and Codes\/}. Springer Universitext, Berlin-Heidelberg, 1993.

\bibitem[5] {}H. Niederreiter and C.P. Xing, {\it Rational Points on Curves over Finite Fields: Theory and Applications\/}. Cambridge University Press, Cambridge, 2001.

\bibitem[6] {}M.A. Tsfasman and S.G. Vladut, {\it Algebraic-Geometric Codes\/}. Dordrecht, The Netherlands: Kluwer, 1991.

\bibitem[7] {}C.P. Xing, {\it Algebraic geometry codes with asymptotic parameters better than the Gilbert-Varshamov and the Tsfasman-Vladut-Zink bounds\/}. IEEE Trans. Inf. Theory, {\bf 47} (2001), 347--352.


\bibitem[8] {}C.P. Xing, {\it Nonlinear codes from algebraic curves improving the Tsfasman-Vladut-Zink bound\/}. IEEE Trans. Inf. Theory, {\bf 49} (2003), 432--437.


\bibitem[9] {}C.P. Xing and H. Chen, {\it Improvements on parameters of one-point AG codes from Hermitian curves\/}. IEEE Trans. Inf. Theory, {\bf 48} (2002), 535--537.

\bibitem[10] {}L. Xu, {\it Improvement on parameters of Goppa geometry codes from maximal curves using the Vladut-Xing method\/}. IEEE Trans. Inf. Theory, {\bf 51} (2005), 2207--2210.

\bibitem[11] {}K. Yang and P.V. Kumar, {\it On the true minimum distances of Hermitian codes\/}. In Lecture Notes in Mathematics. vol. 1518, pp. 99--107, 1991.


\end{thebibliography}
\end{document}